\documentclass[a4paper,12pt]{article}
\usepackage[utf8]{inputenc}
\usepackage{amsthm,amssymb,amsmath}
\newtheorem{theorem}{Theorem}
\newtheorem{lemma}{Lemma}

\newtheorem{corollary}{Corollary}
\newtheorem{question}{Open Problem}

\usepackage{authblk}
\title{Nonrepetitively 3-colorable subdivisions of graphs with a logarithmic number of subdivisions per edge}
\author{Matthieu Rosenfeld\footnote{LIRMM, CNRS, Université de Montpellier}}
\newcommand{\pal}{\rho}
\newcommand{\hpal}{\sigma}
\begin{document}
\maketitle
\begin{abstract}
We show that for every graph $G$ and every graph $H$ obtained by subdividing each edge of $G$ at least $O(\log |V(G)|)$, $H$ is nonrepetitively 3-colorable. In fact, we show that  $O(\log \pi'(G))$ subdivisions per edge are enough, where $\pi'(G)$ is the nonrepetitive chromatic index of $G$.
This answers a question of Wood and improves a similar result of  Pezarski and Zmarz that stated the existence of at least one 3-colorable division with a linear number of subdivision vertices per edge.
\end{abstract}

\section{Introduction}

A sequence $s_1\ldots s_{2n}$ is \emph{a square} if
$s_i=s_{i+n}$  for each $i\in\{1,\ldots,n\}$.
A sequence is \emph{repetitive} if it contains a consecutive subsequence that is a square 
and it is \emph{nonrepetitive} (or \emph{square-free}) otherwise. 
For instance, the words \textbf{hotshots}, \textbf{repetitive} and \textbf{alfalfa} are repetitive and
the words \textbf{total} and \textbf{minimize} are nonrepetitive.

The work of Thue on nonrepetitive words is regarded as the starting point of combinatorics on words
\cite{Thue06,Thue1} 
(see \cite{BerstelThue} for a translation in modern mathematical English).  He showed that there are infinite square-free  sequences over three elements.
Many generalizations and variations of this notion have been studied.
In particular, the notion of nonrepetitive coloring of graphs was introduced
by Alon et al. \cite{alongraph} 
(see \cite{woodSurvey} for a recent survey on this topic).
We say that a coloring (either of the vertices or of the edges) of a graph is \emph{nonrepetitive} if the sequence of colors of any path is nonrepetitive.
The \emph{nonrepetitive chromatic number} (resp. \emph{nonrepetitive chromatic index}) of a graph, denoted by $\pi(G)$ (resp. $\pi'(G)$) is the smallest number of colors in a nonrepetitive coloring of the vertices (resp. the edges) of the graph. 
Alon et al. showed that $\pi'(G)$ is in $O(\Delta^2)$ where $\Delta$ is the maximum degree of $G$ \cite{alongraph}.  
Different authors successively improved the upper bounds on the nonrepetitive chromatic number and the nonrepetitive chromatic index and the best known bound for the nonrepetitive chromatic number is also in $O(\Delta^2)$ \cite{ Dujmovic2016, MontassierEntropie, HARANT2012374, rosenfeldCounting}. 

Nonrepetitive coloring of subdivisions of graphs were also widely studied. 
We say that a graph $G'$ is a \emph{subdivision} of the graph $G$ if $G'$ is obtained by replacing each edge $vw$ of $G$ by a path $P$ with endpoints $vw$, where the new paths are pairwise internally disjoints.
If each edge is replaced by a path with at least $d$ internal vertices then $G'$ is a \emph{$(\ge d)$-subdivision} of $G$.
 Barát and Wood proved that every graph has a nonrepetitively $4$-colorable subdivision \cite{3subdivision4colors}. 
 Pezarski and Zmarz reduced $4$ to $3$ \cite{3subdivision} (solving a conjecture of Grytczuk \cite{Grytczuksubdivisionsurvey}). This result is a strong generalization of Thue's result. 
 In these two results, the number of division vertices per edge is $O(|V(G)|)$ or $O(|E(G)|)$.
 Djumovi\'c et al.showed that every graph has a nonrepetitively $5$-colorable subdivision with $O(\log|V(G)|)$ division vertices per edge \cite{Dujmovic2016}. Their result is in fact stronger than that since it holds in the list-coloring setting and that their bound is in fact $O(\log \Delta)$ where $\Delta$ is the maximal degree of the graph. 
 Finally, Wood proved that every graph has a nonrepetitively $5$-colorable subdivision with $O(\log \pi(G))$ division vertices per edge \cite{woodSurvey}.
It is slightly stronger since it implies the same bound of $O(\log \Delta)$, but it does not hold in the list-coloring setting and requires all the edges to be subdivided in the same amount of internal vertices. 
In a recent survey Wood asked the following question.
 \begin{question}{\cite[Open Question 6.21]{woodSurvey}}
   Does every graph $G$ have a nonrepetitively $3$-colourable subdivision
 with $O(\log|V(G)|))$ or even $O(\log\pi(G))$ division vertices per edge?
 \end{question}
 
 In this article, we give a positive answer to the first part of this question. In fact, we show that there exists a function $f(n)=O(\log n)$ such that any subdivision of any graph $G$ with a least $f(\pi'(G))$ subdivision per edge is nonrepetitively $3$-colorable. 
Since $\pi'(G)$ is in $O(\Delta^2)$, the quantity $O(\log \pi'(G)))$ is smaller than $O(\log \Delta))$ which is itself smaller than the suggested $O(\log|V(G)|))$. However, this is not clear how $O(\log \pi'(G)))$ compares to $O(\log\pi(G))$ so we are not able to solve the second part of this question. The number of subdivision vertices per edge that we require is  $558 \log_2(n)+ O(1)$. 
This result is optimal in the sense that $\Omega(\log n)$ division vertices are needed  on some edges of any nonrepetitively $O(1)$-colorable subdivision of $K_n$\cite{nesetrilsubdivision}. However, we expect the optimal multiplicative coefficient in front of the $\log$ to be much smaller than $558$.
 
Moreover, we show that any subdivision is nonrepetitively colorable as long as each edge is subdivided enough. This is much stronger than showing that there exists one nonrepetitively colorable subdivision.  Remark that all the aforementioned results only showed the existence of one nonrepetitively colorable subdivision. Our proof becomes much simpler if we only care about the existence of one nonrepetitively $3$-colorable subdivision. 
 
The article is organized as follows. We first provide in Section \ref{secpr} some definitions and recall some useful results from the literature. Then in Section \ref{seccow}, we show the existence of a set of words that is needed for our main construction. In particular, we show that there are sets of $n$-good words of exponential size. In Section \ref{secgraph}, we use these sets of $n$-good words to show our main result. The main idea is to start from a nonrepetitive coloring of the edges of a graph and to ``encode'' each color by a well chosen $n$-good word. In the last section, we discuss possible ways to improve the bound on the number of needed subdivisions per edge.

The main idea is implicitly to generalize the notion of square-free morphism.
A \emph{morphism} is a map $h:\Sigma^*\rightarrow\Sigma^*$ such that the image of each word is given by the concatenation of the image of the letters.
A \emph{square-free morphism} is a morphism such that the image of every square-free word is a square-free word. 
A directed graph is nonrepetitively colored if the sequence of colors of every directed path is nonrepetitive.
Given a directed graph $D$, two sets of colors $C$ and $\Sigma$, a nonrepetitive edge coloring $\phi:A(D)\rightarrow C$ of $D$ and a square-free morphism $h:C\mapsto \Sigma$. 
If we can subdivide $D$ in such a way that the sequence of colors of the subdivision of any edge $e$ is $h(\phi(e))$, then this subdivision is nonrepetitively $|\Sigma|$-colorable\footnote{We might need one extra color to color the vertices from the original graph.}.
We cannot use the same construction for undirected graphs. But if we can find a square-free morphism such that the image of every square-free word is square-free even if we replace the image of a letter by its mirror at some arbitrary positions, then we can use the same idea.
There are extra technicalities, since we show that it works for any large enough subdivision.
Naturally, the proof mostly relies on combinatorics on words.

\section{Preliminaries}\label{secpr}
A \emph{word} is a finite sequence over a finite set that we call the \emph{alphabet}. A \emph{factor} of a word is a contiguous subsequence of this word, that is,  if there are two words $p$ and $s$ such that $w=pfs$ then $f$ is a factor of $w$. If $p$ (resp. $s$) is empty then $f$ is also a \emph{prefix} (resp.  \emph{a suffix}) of $w$. A prefix (resp. a suffix) of $w$ is \emph{proper} if it is not equal to $w$.
The \emph{length} of a word $u$ is denoted by $|u|$.
A word $u$ \emph{occurs} in $v$ at position $p$, if the factor of length $|u|$ of $v$ that starts at position $p$ is exactly $u$.
The \emph{mirror image} $\overline{w}$ of a word $w$ is the word obtained by reading $w$ from right to left. We let $\overline{w}^1 = \overline{w}$ and $\overline{w}^0 = w$.

We recall the following result from Shur on the number of square-free ternary words.
\begin{theorem}[\cite{shursqf}]\label{growthsqfree}
 For all $n$, let $C_{sq}(n)$ be the number of ternary square-free words. Then 
 $$\limsup_{n\rightarrow \infty} C_{sq}(n)^{\frac{1}{n}}\ge 1,30175907\,.$$
\end{theorem}
Since any factor of a square-free word is square-free, $C_{sq}$ is a submultiplicative function (i.e., we have $C_{sq}(i+j)\le C_{sq}(i)C_{sq}(j)$,  for all $i,j$). 
By Fekete's Lemma we deduce the following Corollary of Theorem \ref{growthsqfree}.
\begin{corollary}\label{explbonsqfree}
  For all  integer $n\ge1$,  $C_{sq}(n)> 1.3^n$. 
\end{corollary}

We now recall Tur\'an's Theorem and a simple corollary.
\begin{theorem}[Tur\'an's Theorem]
Let $G$ be any graph with $n$ vertices, such that $G$ is $K_{r+1}$-free. Then the number of edges in $G$ is at most 
$$|E(G)|\le \left(1-\frac{1}{r}\right)\frac{n^2}{2}\,.$$ 
\end{theorem}
\begin{corollary}\label{exindepset}
Any graph $G$ contains independent set of size at least $\frac{n}{1+d(G)}$ where $d(G)$ is the average degree of $G$.
\end{corollary}
\begin{proof}
 Let $G$ be a graph of average degree $d(G)$ and let $H$ be the complement of $G$. 
 Then the number of edges of $H$ is
\begin{align*}
 |E(H)|=\frac{n(n-1)}{2}-|E(G)|=\frac{n(n-1)}{2} -\frac{nd(G)}{2}&=\frac{n^2}{2}\left(1-\frac{1+d(G)}{n}\right)\\
 &>\frac{n^2}{2}\left(1-\frac{1}{\frac{n}{1+d(G)}-1}\right)\,.
\end{align*}
 By Tur\'an's Theorem, this implies that $H$ cannot be $K_{\frac{n}{1+d(G)}}$-free. Hence, $G$ contains an independent set of size $\frac{n}{1+d(G)}$ as desired.
\end{proof}

\section{Good sets}\label{seccow}
 
Let $\hpal=1202120121021201021$, $\overline \hpal=1201021201210212021$ be the mirror image and let $\pal=\overline \hpal 0 \hpal$. By construction $p$ is a palindrome, $|\hpal|=19$ and $|\pal|=39$.

A word $v$  is \emph{nice} if $\pal v\pal$ is square-free and contains only two occurrences of $\hpal$ and two occurrences of $\overline{\hpal}$.
The only occurrences of $\hpal$ and $\overline \hpal$ are inside the two occurrences of $\pal$. Let $\mathcal{N}$ be the set of nice words and for any integer $n$, let $\mathcal{N}_n$ be the set of nice words of length $n$.

For all $n\ge 8750$, let $l_n$ be the lexicographically least word of the set $\mathcal{N}_n$. We show in Lemma \ref{sizengoodset} that there exist nice words of every length at least $8750$. So $l_n$ is properly defined for any $n\ge 8750$.

A set of words $S$ is \emph{mirror-free} if for any word $w$ from $S$ such that $w\not=\overline w$,  $\overline w$ is not in $S$. For any integer $n\ge 8750$, a set of words $\mathcal{S}\subseteq\mathcal{N}_n$ is \emph{$n$-good} if it is mirror-free and if
for all $i\in\{2n+100,\ldots,7n\}$ and all $u,v\in \mathcal{S}$ such that $u\not=v$  the words $\pal u\pal v\pal l_i$, $\pal\overline u\pal v\pal l_i$, $\pal u\pal \overline v\pal \overline{l_i}$ and $\pal\overline u\pal\overline v\pal\overline{l_i}$ are square-free.

In Lemma \ref{mainlemmainject}, we show the central property of $n$-good set, but first we show the following Lemma as a warm-up exercise.

\begin{lemma}\label{warmup}
Let $n\ge 8750$ be an integer, $S$ be an $n$-good set, $u,v\in S$ and $i\in \{2n+100,\ldots, 7n\}$ be an integer then $\pal u\pal l_i\pal v\pal$ is square-free.
\end{lemma}
\begin{proof}
Suppose, for the sake of contradiction, that $\pal u\pal l_i\pal v\pal$ contains a square $ww$.
First remark, that since $S$ is an $n$-good set $\pal u\pal l_i\pal$ and
$\pal l_i\pal v\pal$ are both square-free.
Hence there exist a non-empty suffix $x$ of $\pal u$ and a non-empty prefix $y$ of $v\pal$ such that $x\pal l_i\pal y$ is a square.

Now $l_i\pal$ is not a factor of $w$ since it contains only one occurrence of $\pal$ and that $l_i$ is much longer than the gap between any other occurrences of $\pal$. For the same reason $\pal l_i$ is not a factor of $w$.
There exists $x',y'$ such that $l_i=x'y'$ and $w=x\pal x'=y'\pal y$.
Since $|x'|+|y'|=|l_i|\ge 2n+100$, assume, without loss of generality that $|x'|\ge n+50$.
Then $|x'|> | y|$ and $xp$ is a proper prefix of $y'\pal$ which implies that there is a second occurrence of $\pal$ in $y'\pal$. This contradicts the fact that $l_i$ is a nice word. 
\end{proof}

\begin{lemma}\label{mainlemmainject}
Let $k$ and $n$ be two integers. Let $S$ be an $n$-good set, $\Sigma$ be an alphabet and $f:\Sigma\mapsto S$ be an injective map. Then, for any square-free word $w_1\ldots w_k\in \Sigma^k$, any sequence of integers $(s_i)_{1\le i\le k}\in \{2n+100,\ldots, 7n\}^k$ and any sequence $(r_i)_{1\le i\le k}\in \{0,1\}^k$ the word 
$$\pal\prod_{i=1}^k \overline {f(w_i)\pal l_{s_i}\pal f(w_i)}^{r_i}\pal=
\pal\overline {f(w_1)\pal l_{s_1}\pal f(w_1)}^{r_1}
\ldots 
\pal\overline {f(w_k)\pal l_{s_k}\pal f(w_k)}^{r_k}\pal$$
is square-free.
\end{lemma}
\begin{proof}
Suppose, for the sake of contradiction, that there is a square $uu$ in $\pal\prod_{i=1}^k \overline {f(w_i)\pal l_{t_i}\pal f(w_i)}^{r_i}\pal$.
There are $x_1,\ldots,x_{3k}\in\mathcal{N}_n$ such that
for all $i$, $\overline {f(w_i)\pal l_{s_i}\pal f(w_i)}^{r_i}=x_{3i-2}\pal x_{3i-1}px_{3i}$. This implies that 
$$\pal\prod_{i=1}^k \overline {f(w_i)\pal l_{s_i}\pal f(w_i)}^{r_i}\pal=\pal x_1\pal x_2\pal \ldots \pal  x_{3k}\pal $$
and that there are no other occurrences of $\hpal $ and $\overline{\hpal} $ than the $3k+1$ occurrences that are inside the occurrences of $\pal $.

Let  $L= \left\{l_i:i\in \left\{2n+100,\ldots, 7n\right\}\right\}$. Then for any $j$,  $x_j\in L$ if and only if $j\equiv 2 \mod 3$.
Since $S$ is mirror-free, for any $i$, $x_i$ and $x_{i+1}$ are different. 
Thus, by the definition of $n$-good sets and by Lemma \ref{warmup}, for any $i$, $\pal x_i\pal x_{i+1}\pal x_{i+2}\pal$ is square-free. Thus there are at least $3$ occurrences of $\pal$ in $uu$. At least one of $\hpal$ or $\overline{\hpal}$ appears twice in $u$.
Assume, without loss of generality, that the middle of the square $uu$ does not cut any occurrence of $\hpal $ (if it is not the case we can use the same argument with $\overline{\hpal} $ instead). Then $\hpal$ occurs at least twice in $u$.

Let $l\ge2$ be the number of occurrences of $\hpal $ in $u$. 
Since the explicit occurrences of $\hpal$ are the only occurrences, and that there are as many occurences of $\hpal$ in the two consecutives occurences of $u$,
there exist an integer $i$ and words $y,y',z, z'$ such that $y$ is a proper suffix of $\hpal x_{i-l}\overline{\hpal} 0$, $y'z=x_i\overline{\hpal}0$, $z'$ is a prefix of $x_{i+l}\overline{\hpal}  0 $ and  
\begin{align*}
u&=y\hpal x_{i+1-l}\overline{\hpal}  0 \hpal x_{i+2-l}\overline{\hpal}  0 \hpal \ldots \overline{\hpal}  0 \hpal x_{i-1}\overline{\hpal}  0 \hpal y'\\
&=z\hpal x_{i+1}\overline{\hpal}  0 \hpal x_{i+2}\overline{\hpal}  0 \hpal \ldots\overline{\hpal}  0 \hpal  x_{i+l-1}\overline{\hpal}  0 \hpal z'\,.                                                                                                             \end{align*}
Moreover, the occurences of $\hpal$ have to match each others and \emph{synchronise} the rest of $u$, that is, $y=z$, $y'=z'$ and for all $j\in\{0,\ldots,l-2\}$, 
\begin{equation}\label{syncrhonisedxi}
x_{i+1-l+j}=x_{i+1+j}\,.
\end{equation}

Recall that, for any $j$,  $x_j\in L$ if and only if $j\equiv 2 \mod 3$.

If $l=2$, then $x_{i-1}=x_{i+1}$ which implies $i\equiv 2\mod3$.
Moreover,  $y=z$ and $y'=z'$ imply 
$$|x_i\overline{\hpal} 0|=|y|+|z'|< |\hpal x_{i-l}\overline{\hpal} 0|+|x_{i+l}\overline{\hpal}  0|= |x_{i-l}|+|x_{i+l}|+59\,.$$ But since, $i\equiv 2\mod3$, $|x_i|\ge 2n+100$ and $|x_{i+1}|=|x_{i-1}|=n$ which contradicts the previous equation.

Let us now take care of the case $l\ge3$.
We first show that $l\equiv 0 \mod 3$.
For all  $j\in\{0,\ldots,l-2\}$,   $i+1-l+j\equiv 2 \mod 3$ if and only if $i+1+j\equiv 2 \mod 3$.
Thus if  $i+1-l\equiv 2 \mod 3$ or $i+2-l\equiv 2 \mod 3$ then $l$ has to be divisible by $3$. 
If it is not the case then  $i+1-l\equiv 0 \mod 3$, $|x_{i+1-l}|=|x_{i+2-l}|=n$ which implies $|x_{i+1}|=|x_{i+2}|=n$ and $i+1\equiv 0 \mod 3$ and finally $l \equiv 0\mod 3$. 

We have three different cases to consider.
\paragraph{Case $i+1-l \equiv 1\mod 3$:} First recall that, for every $j\equiv 1\mod 3$, $x_j= \overline {f(w_{(j+2)/3})}^{(j+2)/3}$. Thus for all integer $j$, 
$x_{i+1+3j-l}= \overline {f(w_{(i-l)/3+1+j})}^{r_{(i-l)/3+1+j}}$ and 
$x_{i+1+3j}=\overline { f(w_{i/3+1+j})}^{r_{i/3+1+j}}$.
By equation \eqref{syncrhonisedxi}, for all $j\in\{0,\ldots, l/3-1\}$,
 $$\overline {f(w_{(i-l)/3+1+j})}^{r_{(i-l)/3+1+j}}=\overline { f(w_{i/3+1+j})}^{r_{i/3+1+j}}\,.$$
The function $f$  is an injective and maps to $S$ which is mirror-free, so the previous equation implies that for all  $j\in\{0,\ldots, l/3-1\}$, $$w_{(i-l)/3+1+j}=w_{i/3+1+j}\,.$$
This implies that there is a square in $w$ which is a contradiction.

\paragraph{Case $i+1-l\equiv 0 \mod 3$:} We use the same idea as in the previous case with the fact that for every $j\equiv 0\mod 3$, $x_j= \overline {f(w_{j/3})}^{r_{j/3}}$. In this case we obtain that for all  $j\in\{0,\ldots, l/3-1\}$, $$w_{(i+1-l)/3+j}=w_{(i+1)/3+j}\,.$$
This implies that there is a square in $w$ which is a contradiction.

\paragraph{Case $i+1-l\equiv 2 \mod 3$:} This case is almost identical to the previous ones. We know that for every $j\equiv 0\mod 3$, $x_j= \overline {f(w_{j/3})}^{r_{j/3}}$. Moreover for all  $j$,
$i+2-l+3t\equiv 0 \mod 3$. By the same argument, for all  $j\in\{0,\ldots, l/3-1\}$, $$w_{(i+1-l)/3+j}=w_{(i+1)/3+j}\,.$$
This implies that there is a square in $w$ which is a contradiction.
\end{proof}

This property is essential to construct the nonrepetitive coloring of a subdivided graph. The idea is to encode the colors of the edges of a nonrepetitive edge coloring of the initial graph. Any vertex from the initial graph will be colored by $0$ and the path corresponding to any edge colored $c$ in the original graph should receive the color sequence $\hpal f(c) \pal l_i \pal f(c) \overline{\hpal}$ (with the $l_i$ of the right length). The fact that we can choose a different $l_i$ for every edge means that we can find a right encoding as long as the edge is subdivided enough. 
The fact that we can replace the encoding of each $w_i$ by the mirror image of the encoding means that we can take an arbitrary orientation of the subdivided edge to apply the encoding to the corresponding path.

We also need a variant of this property. This variant will be useful for paths that start and end in the subdivision of the same edge (i.e., paths that appear in the subdivision of a cycle).

\begin{lemma}\label{mainlemmainjectbis}
Let $k$ and $n$ be two integers. Let $S$ be an $n$-good set, $\Sigma$ be an alphabet and $f:\Sigma\mapsto S$ be an injective map.
Let  $w_1\ldots w_k\in \Sigma^k$ be a square-free word such that $w_2w_3\ldots w_kw_1$ is also square-free. Let $(t_i)_{1\le i\le k}\in \{2n+100,\ldots, 7n\}^k$ be a sequence of integers, and  $(r_i)_{1\le i\le k}\in \{0,1\}^k$ be a sequence of $0$ and $1$.
Let  $a$ and $b$ be a pair of words such that $\overline {f(w_1)\pal l_{t_1}\pal f(w_1)}^{r_1}\pal = ab$. Then the word
\begin{align*}
&b\left(\prod_{i=2}^k \overline {f(w_i)\pal l_{t_i}\pal f(w_i)}^{r_i}\pal \right) a\\
&=
b \pal \overline {f(w_2)\pal l_{t_2}\pal f(w_2)}^{r_2}\pal 
\overline {f(w_2)\pal l_{t_2}\pal f(w_2)}^{r_2}\pal 
\ldots 
\pal \overline {f(w_k)\pal l_{t_k}\pal f(w_k)}^{r_k}\pal  a
\end{align*}
is square-free.
\end{lemma}
\begin{proof}
  Suppose, for the sake of contradiction, that there is a square $uu$ in $b\left(\prod_{i=2}^k \overline {f(w_i)\pal l_{t_i}\pal f(w_i)}^{r_i}\pal \right) a$.
  
  Since $w_1w_2\ldots w_k$ and  $w_2w_3\ldots w_kw_1$ are both square-free, Lemma \ref{mainlemmainject} implies that
  $b\left(\prod_{i=2}^k \overline {f(w_i)\pal l_{t_i}\pal f(w_i)}^{r_i}\pal \right)$ and $\left(\prod_{i=2}^k \overline {f(w_i)\pal l_{t_i}\pal f(w_i)}^{r_i}\pal \right)a$ are also square-free. Hence the first occurrence of $u$ starts in $b$ and the second occurrence ends in $a$. Assume, without loss of generality, that the middle of the square does not cut an occurrence of $\hpal $.
  Since $k>3$, we know that $uu$ contains enough occurrences of $\hpal $ to be synchronized by the occurrences of $\hpal $ and this implies that the number of occurrences of $\hpal $ in $u$ is a multiple of $3$. With the same argument as in the proof of Lemma \ref{mainlemmainject} one easily shows that for every $i\in\left\{1,\ldots,\frac{k}{2}-1\right\}$,
  $$w_{i+1}= w_{k/2+i+1}\,.$$
  We also easily verify that  $\overline{f(w_{k/2+1})\pal l_{t_{k/2+1}}\pal f(w_{k/2+1})}^{r_{k/2+1}}$ is a factor of $ab=\overline {f(w_1)\pal l_{t_1}\pal f(w_1)}^{r_1}\pal $. This is only possible if $w_1 = w_{k/2+1}$. Thus, for every $i\in\left\{0,\ldots,\frac{k}{2}-1\right\}$,
  $$w_{i+1}= w_{k/2+i+1}\,.$$
  This is a contradiction since $w_1\ldots w_k$ is square-free.  
\end{proof}

\subsection{Exponentially large good sets}
We show in this subsection that there are exponentially many nice words of any length and we use that to show that there are exponentially large $n$ good sets.

Let $h:\Sigma\mapsto\Sigma^*$, be the map such that
\begin{align*}
h(0)=&\{
        012102120210201 021201210,
        0121021202102012021201210\}\,,\\
h(1)=& \{120210201021012102012021, 1202102010210120102012021\},\\
h(2)=& \{201021012102120210120102, 2010210121021201210120102\}\,.
\end{align*}
A word $v$ is an \emph{image} of $w$ by $h$ if $v$ can be obtained by replacing each occurrence of any letter $i$ of $w$ by any word of the corresponding set $h(i)$. The set of images of $w$ by $h$ is denoted by $h(w)$.
The authors of  \cite{nicesubstitution} introduced $h$ and showed the following property.\footnote{In fact, the map that they used contains $4$ words in each set and here we used only 2 of them.}
\begin{theorem}[{\cite[Theorem 20]{nicesubstitution}}]\label{sqfreesub}
For any square-free word $w\in\{0,1,2\}^*$ and any $v\in h(w)$,
$v$ is square-free.
\end{theorem}

We will use $h$ to show that the set of nice words has exponential growth.
First, we need  a few simple facts about $h$, $\hpal$ and $\pal$.

\begin{lemma}\label{observations}
 \begin{enumerate}
 \item For any letters $a,b,c\in\{0,1,2\}$, $v\in h(ab)$ and $v'\in h(c)$, the word $v'$ cannot occur as an internal factor of $v$ (i.e., the only possible occurrences of $v'$ are as suffix or prefix).
  \item For any $a,b\in\{0,1,2\}$  and $v\in h(ab)$,
 neither  $\hpal$ or $\overline{\hpal}$  are factor of $v$.
 \item If $a\in\{1,2\}$, $b\in \{0,1,2\}\setminus\{a\}$ and $v\in h(0ab)$, then the word $\pal v$ is square-free and contains exactly one occurrence of $\hpal$ and $\overline{\hpal}$.
 \item If $a\in\{1,2\}$, $b\in \{0,1,2\}\setminus\{a\}$ and $v\in h(ba0)$, then the word $v\pal$ is square-free and contains exactly one occurrence of $\hpal$ and $\overline{\hpal}$.
 \end{enumerate}
\end{lemma}

 It is a bit tedious to verify the four claims of this Lemma by hand so we provide a simple computer program that verifies Lemma \ref{observations}\footnote{See the ancillary file  \texttt{verifying\_lemma4.cpp}.}.
The first of this 4 facts implies that images by $h$ synchronize, that is, as long as a factor of an image is of length at least $50$ (the length of the images of two letters), there is a unique way to split it into different images by $h$. This kind of properties are really useful to establish that an image of a square-free word by $h$ does not contains any large squares (i.e., squares that are large enough to allow us to use this synchronization property). On the other hand, the other facts are quit useful to establish that there are no short squares. With this lemma in hand it is relatively simple to provide a proof of Theorem \ref{sqfreesub} (which we will not do). 
We can use these facts to show that there are exponentially many nice words.

\begin{lemma}\label{hissqfreeandpfree}
Let $w\in \{0,1,2\}^*$ be a word such that $|w|\ge4$ and  $10w01$ is square-free. Then any $v\in h(0w0)$ is nice.
\end{lemma}
\begin{proof}
Let $w$ and $v$ be as in the Theorem statement.
Let $w_1\ldots w_n=w$ with $w_1,\ldots,w_n\in \{0,1,2\}$ and for all $i$ let $v_i\in h(w_i)$ such that $v=v_1\ldots v_n$.
Recall that we need to show that $\pal v\pal$ is square-free and contains only two occurrences of $\hpal$ and two occurrences of $\overline{\hpal}$. 
By 2., 3.  and 4. of Lemma \ref{observations}, any occurrence of $\hpal$ or $\overline{\hpal}$ is in $\pal$, so it only remains to show that $\pal v\pal$ is square-free.
 
 Suppose, for the sake of contradiction, that there is a square in $\pal v\pal$, that is, there are words $x,y\in\{0,1,2\}^*$ and $u\in\{0,1,2\}^+$ such that $\pal v\pal =xuuy$. Theorem \ref{sqfreesub} implies that $uu$ cannot be a factor of $v$. One easily verifies that $\pal$ is square-free. 
 Thus $x$ is a proper prefix of $\pal$ or $y$ is a proper suffix of $\pal$. Assume, without loss of generality, that $x$ is a proper suffix of $\pal$.
 Let $r$ be the nonempty suffix of $\pal$ such that $xr=\pal$.
 Fact 3. of Lemma \ref{observations} implies that the square $uu$ is not a factor of $\pal v_1v_2v_3$, thus $|uu|\ge |v_1v_2v_3|+2$ and $v_1$ is a factor of $u$.
 Thus $rv_1$ is a prefix of $u$.  
 We have to distinguish between two different cases depending on the length of $r$.

 \paragraph{Case $|r|\ge4$:}
By hypothesis $|w|\ge4$ and $|v_2\ldots v_n|\ge 5\times24>|v_1|+2|\pal|\ge |v_1r\pal |$ which gives $|rv_1|+\frac{|r+v+\pal|}{2}<|rv|$.
We deduce that the $|rv_1|$ first letters of the second occurrence of $u$ do not overlap with the final occurrence of $\pal$. 
 This implies that $rv_1$ is a factor of $v$. 
 From Fact 1. of Lemma \ref{observations}, $v_1$ can only appear as the image of $0$ and thus $r$ must appear as the suffix of the image of a letter. 
 However, it is easy to verify that $r$ is not the suffix of any image of a letter if $|r|\ge4$ (it is enough to verify this with $|r|=4$), which is a contradiction.

 \paragraph{Case $|r|\le3$:}
 Then $r$ is also the suffix of an image of $1$ by $h$ (since $021$ is suffix of any image of $1$). 
There is $v_0\in h(1)$ such that $uu$ is also a factor of $v_0v\pal$. By hypothesis $w$ was chosen such that $10w01$ is square-free and Theorem \ref{sqfreesub} implies that there is no square in $v_0v$. Thus the square in $v_0v\pal$ overlaps with the final occurrence of $\pal$. By symmetry of the previous case, the square overlaps by at most $3$ letters with the final occurence of $\pal$ which implies that the square is also an image of $10w01$ which is a contradiction since any image of $10w01$ by $h$ is square-free by Theorem \ref{sqfreesub}.
\end{proof}
We can deduce an exponential lower bound on the size of $\mathcal{N}_n$. 

\begin{theorem}\label{expboundsqpfree}
 For all $n>8750$,
 $$|\mathcal{N}_n|> \frac{1.01^{n}}{7.8}\,.$$
\end{theorem}
\begin{proof}
Let $P$ be the set of words over $\{0,1,2\}$ such that for any $w\in P$,  $|w|\ge4$, $10w01$ is square-free.
Lemma \ref{hissqfreeandpfree} implies that for all $n$
 $$|\mathcal{N}_n|\ge \left|\left\{v\in \{0,1,2\}^n | \exists w \in P, v\in h(0w0) \right\}\right|\,.$$
The set of images by $h$ of two different words are distinct, hence
\begin{equation}
|\mathcal{N}_n|\ge \left|\left\{ w \in P|h(0w0)\cap\{0,1,2\}^n\not=\emptyset \right\}\right|\,.\label{eqboundSn}
\end{equation}

Every letter has an image of length $24$ and an image of length $25$ by $h$. Thus for any integer $n\ge 14\times 25^2=8750$ and any word $w$ of length $\frac{n}{24}-13\le|w|\le\frac{n}{24}$, $w$ admits at least an image of size $n$ by $h$.
That is, for any  $n\ge 8750$,
$$\{ w \in P|h(0w0)\cap\{0,1,2\}^n\not=\emptyset \ \}|\ge
\left|\left\{w\in P|  \frac{n}{24}-13\le|w|\le\frac{n}{24}\right\}\right|\,.$$

By symmetry, there are exactly $\frac{C_{sq}(n)}{6}$ square-free words over $\{0,1,2\}$  of length $n$ starting by $10$.
Moreover, one easily verifies that every square-free word over $\{0,1,2\}$ of length $14$ contains at least one occurrence of $01$.
Thus for every integer $n$, there are at least $\frac{C_{sq}(n)}{6}$ ternary square-free words of length between $n-13$ and $n$  starting with $10$ and ending with $01$.

We can now apply Corollary \ref{explbonsqfree},
 $$|\{ w \in P|h(0w0)\cap\{0,1,2\}^n\not=\emptyset \ \}|\ge\frac{C_{sq}(\lfloor\frac{n}{24}\rfloor)}{6}> \frac{1.3^{\lfloor\frac{n}{24}\rfloor}}{6}> \frac{1.01^{n}}{7.8}\,.$$
 Together with equation \eqref{eqboundSn}, we conclude 
 $|\mathcal{N}_n|> \frac{1.01^{n}}{7.8}$.
\end{proof}

We use the exponential lower bound to establish the existence of exponentially large $n$-good sets, but first we show one more property of nice words.

\begin{lemma}\label{shortsquareinpref}
Let $n$ be a positive integer and $u,v\in \mathcal{N}$. If the word $\pal u\pal v\pal$ is not square-free then $u$ is a prefix of $v$ or $v$ is a suffix of $u$. 
In particular, if $|u|=|v|$ then $u=v$.
\end{lemma}
\begin{proof}
Let $n$ be a positive integer and $u,v\in \mathcal{N}_n$ such that $u\not=v$.
Let $ww$ be a square in $\pal u\pal v\pal$.
Since $u$ and $v$ are nice, $\pal u\pal$ and $\pal v\pal$ are square-free. Hence the second occurrence of $\pal$ is a factor of $ww$. We also know that the only occurrences  of $\hpal$ in $\pal u\pal v\pal$ (resp. of $\overline{\hpal}$) are the three occurrences inside each occurrence of $\pal$.

Suppose, for the sake of contradiction, that there are two non-empty words $u_1$ and $u_2$ such that $u=u_1u_2$, $w$ is a suffix of $\pal u_1$ and a prefix of $u_2\pal v\pal$. Since $w$ contains $\pal$ as a factors and that there are exactly three occurrences of $\pal$ this implies that $w= \pal u_1$ which is a contradiction with the fact that $u_2\pal v\pal$  does not starts with $\pal$. By symmetry, we reach a similar contradiction if we try to split $v$ in $v_1v_2$.

Hence there exists $\pal_1$, $\pal_2$ such that $\pal_1\pal_2=\pal= \overline{\hpal} 0 \hpal$ and $w$ is a suffix of $\pal u\pal_1$ and a prefix of $\pal_2v\pal$.
Assume, without loss of generality, that $\hpal$ is a suffix of $\pal_2$ (otherwise $\overline{\hpal}$  is a prefix of $\pal_1$ and the rest of the argument is symmetric).
Since the only occurrence of $\hpal$ in $\pal u\pal_1$ is inside $\pal$, we deduce that $\pal_2 u\pal_1=w$. Since $w$ is a prefix of $\pal_2v\pal$, $u$ is a prefix of $v$ or $v$ is a prefix of $u$.
\end{proof}

We can finally show the existence of exponentially large $n$-good sets.
\begin{lemma}\label{sizengoodset}
 For any $n>8750$, there exists an $n$-good set $S$ of size at least
  $$|S|\ge\frac{1.01^{n}}{16(60n^2+1)}\,.$$
\end{lemma}
\begin{proof}
Let $\mathcal{N}'_n$ be the set obtained by removing from $\mathcal{N}_n$ any prefix or suffix of
every $l_i$ with  $i\in\{2n+100,\ldots,7n\}$ and by keeping for each pair of mirror images only the lexicographically smallest of the two.
Each $l_i$ is responsible for removing at most two words from $\mathcal{N}_n$ so $|\mathcal{N}'_n|\ge \frac{|\mathcal{N}_n|-10n}{2}\ge\frac{1.01^{n}}{15.6}-5n$. Since $n>8750$, we can simplify the bound $$|\mathcal{N}'_n|\ge\frac{1.01^{n}}{16}\,.$$

For any $u,v\in\mathcal{N}'_n$, we say that $u$ forbids $v$ if $u\not=v$ and for some $i\in\{n+1,n+2,\ldots,5n\}$, $\pal u\pal v\pal l_i\pal$, $\pal\overline{u}\pal v\pal l_i$, $\pal u\pal \overline{ v}\pal \overline{l_i}\pal$ or $\pal\overline{u}\pal\overline{v}\pal\overline{l_i}\pal$ contains a square.

We now count how many words $v$ are forbidden by a given $u$.

Let $u,v\in\mathcal{N}'_n$ be such that $u\not=v$ and  $\pal u\pal v\pal l_i\pal$ is not square-free. 
Lemma \ref{shortsquareinpref} implies that both $\pal u\pal v\pal$ and $\pal v\pal l_i\pal$ are square-free (since $u\not=v$ and $v$ is not a prefix of $l_i$).
There is a non-empty suffix $u'$ of $\pal u$ and a non-empty prefix $l'$ of $l_i\pal$ such that the square $ww= u'\pal v\pal l'$. Moreover, there exist two non-empty words $v_1$ and $v_2$ such that $\hpal v\overline{\hpal}=v_1v_2$ and $w=u'\overline{\hpal} 0 v_1=v_2 0\hpal l'$. Indeed, the middle of the square cannot be located outside of $\hpal v\overline{\hpal}$ since there would be too many occurences of $\hpal$ or $\overline{\hpal}$ on one side of the square.
Finally, remark that either $v_1$ contains $\hpal$ as a prefix or $v_2$ contains $\overline{\hpal}$ as a suffix (both could be true). In both cases, using the fact that there are only two other occurrences of $\hpal$ and $\overline{\hpal}$, we deduce that $|w|=|v\pal |$. Thus $v_1$ is a prefix of $0\hpal l_i\pal$ and $v_2$ is a suffix of $\pal u\overline{\hpal} 0$ and $v$ is uniquely determined by $u$, $l_i$ and the position of the square.
There are $|v\pal |=n+39$ possible positions, less than $5n$ possibles values for $l_i$, so $u$ forbids at most $(n+39)\times 5n$ words because of $\pal u\pal v\pal l_i$. The count is similar for $\pal\overline{u}\pal v\pal l_i$, $\pal u\pal \overline{v}\pal \overline{l_i}\pal$ and $\pal\overline{u}\pal\overline{v}\pal\overline{l_i}\pal$, so  $u$ forbids at most $(n+39)\times 5n\times 4$ words. This is upper bounded by $30n^2$ since $n>8750$.

Let $G$ be the graph whose vertices are the words from $\mathcal{N}'_n$ and such that two words share an edge if one of them forbids the other one. 
The set of words corresponding to any independent set of $G$ is an $n$-good set.
Let $S$ be the set of words corrsponding to the largest independent set of $G$.
Since every word forbids at most $30n^2$ words, the average degree of the vertices of $G$ is at most $60n^2$.
By Corollary \ref{exindepset}, there is an independent set of size at least $\frac{|\mathcal{N}'_n|}{60n^2+1}$.
Thus $|S|\ge \frac{1.01^{n}}{16(60n^2+1)}$.
\end{proof}

\section{The final construction}\label{secgraph}
A graph $G'$ is a \emph{$(\ge a,\le b)$-subdivision} of a graph $G$ if $G'$ can be obtained by subdividing each edge of $G$ in at least $a$ and at most $b$ division vertices.

In Lemma \ref{lemmaConstruction}, we use our results on $n$-goods sets to show that, if each edge of the graph is subdivided enough, but not too much, then we can nonrepetitively $3$-color the resulting graph. To obtain Theorem \ref{mainth}, we then show that we can easily handle the edges that have too many subdivision vertices.

\begin{lemma} \label{lemmaConstruction}
  Let $G$ be a graph and $n\ge 8750$ an integer such that $\pi'(G)\le \frac{1.01^{n}}{16(60n^2+1)}$.
  Then for any $(\ge 4n+216, \le 9n)$-subdivision $G'$ of $G$, 
  $$\pi(G')=3\,.$$  
\end{lemma}
\begin{proof}
Let $n$, $G$ and $G'$ be as in the lemma statement.
Let $C$ be a set of colors of size  $\pi'(G)$ and $\phi$ be a nonrepetitive edge $C$-coloring of $G$.

By Lemma \ref{sizengoodset}, there is an $n$-good set $S$ such that $|S|\ge\pi'(G)$.
Let $f$ be an injective map from $C$ to $S$. Let $\vec o$ be an arbitrary orientation of the edges of $G$.

Let $\phi':V(G')\mapsto \{0,1,2\}$ be the 3-coloring of the vertices of $G'$ such that 
\begin{itemize}
    \item the color of every vertex of $G'$ that corresponds to an original vertex of $G$ has color $0$,
    \item for any edge $e$ from $G$ subdivided in $(v_1,\ldots,v_k)$ in $G'$ with the $v_i$ ordered according to $\vec o(e)$, the sequence of colors $(\phi'(v_1),\ldots,\phi'(v_n))$ is equal to $\hpal f(\phi(e))\pal l_{k-116-2n}\pal f(\phi(e))\overline{\hpal}$.
\end{itemize}
Remark that $|\hpal f(\phi(e))\pal|+ |pf(\phi(e))\overline{\hpal} |= 2|\pal|+2|\hpal|+2n= 2n+116$ and thus 
$|\hpal f(\phi(e))\pal l_{k-116-2n}\pal f(\phi(e))\overline{\hpal}| =k$. Thus $\phi'$ is well-defined. Our goal is now to show that $\phi'$  is nonrepetitive.

First remark that since every edge of $G$ is subdivided at least $4n+216$ times and at most $9n$ times this implies that for each edge of $G$ subdivided into $k$ vertices, $2n+100\le k-116-2n\le 7n$. So the length of the $l_i$ allows us to apply Lemma \ref{mainlemmainject} and Lemma \ref{mainlemmainjectbis}.

Let $\mathbf{p}$ be a path in $G'$ whose two extremities do not belong to the subdivision of the same edge of $G$. Then it is a subpath of the subdivision of some path in $G$.
Let $e_1,\ldots, e_k$ be this path of $G$. For all $i\in\{1,\ldots, k\}$, let $w_i=\phi(e_i)$, let $r_i$ be $0$ if $\vec o (e_i)$ goes in the same direction as the orientation of the path and $r_i=1$ otherwise. For all $i\in\{1,\ldots, k\}$, let $d_i$ be the integer such that $e_i$ is subdivided into $d_i$ vertices in $G'$ and let $t_i= d_i -116-2n$. 
Then by definition the sequence of colors of the path $\mathbf{p}$ from $G'$ is a factor of
$$\pal\prod_{i=1}^k \overline {f(w_i)\pal l_{t_i}\pal f(w_i)}^{r_i}\pal\,.$$
Moreover, since $\phi$ is nonrepetitive, $w_1\ldots w_p$ is square-free. 
By Lemma \ref{mainlemmainject}, $\mathbf{p}$ is nonrepetitively colored by $\phi'$.

Now we need to show that the same property holds if the two extremities of a path $\mathbf{p}$ of $G'$ belong to the subdivision of the same edge. If the path is short and completely contained in an edge then this is in fact solved as the previous case. Then the remaining case is that 
$\mathbf{p}$ starts in the subdivision of an edge $e_1$ of $G$, then leaves this subdivision, and comes back to it by the other side. Let $e_1,e_2,\ldots, e_n, e_1$ be the edges of $G$ whose subdivision contains $\mathbf{p}$.
For all $i\in\{1,\ldots, k\}$, let $w_i=\phi(e_i)$. For all $i\in\{1,\ldots, k\}$, define $r_i$ and $t_i$ as in the previous case. Then there are two words $a$ and $b$ such that $\overline {f(w_1)\pal l_{t_1}\pal f(w_1)}^{r_1}\pal = ab$ and such that the sequence of colors of $\mathbf{p}$ is
$$b\left(\prod_{i=2}^k \overline {f(w_i)\pal l_{t_i}\pal f(w_i)}^{r_i}\pal \right) a\,.$$
By Lemma \ref{mainlemmainjectbis}, $\mathbf{p}$ is nonrepetitively colored by $\phi'$.

We showed that every possible path of $G'$ is nonrepetitively colored by $\phi'$ which implies that $\phi'$ is a nonrepetitive $3$-coloring of $G'$.
\end{proof}

\begin{lemma}\label{lemmasubddoesnotincreasepi}
  Let $G$ be a graph and $H$ be a subdivision of $G$ then
  $\pi'(H) \le 2\pi'(G)+3$.
\end{lemma}
\begin{proof}
  Let $\phi$ be a nonrepetitive edge coloring of $G$ over the set of colors $C$ of size $\pi(G)$.
  Let $C'$ be the set of colors obtained by adding three new colors $\alpha,\beta,\gamma$ and for each color $c\in C$ a new color $c'$.
  
  Let $\phi'$ be an edge coloring of $H$ such that for each edge $e$ of $G$:
  \begin{itemize}
    \item if $e$ is not subdivided in $H$ then it has the same color in $H$ and in $G$,
    \item if $e$ is subdivided into two edges $e_1$ and $e_2$ then $\phi'(e_1)=\phi(e)$ and $\phi'(e_2)=\phi(e)'$,
    \item if $e$ is subdivided in $k\ge3$ edges $e_1,\ldots, e_k$, then $\phi'(e_1)=\phi'(e_k)=\phi(e)$ and the sequence $\phi'(e_2)\ldots\phi'(e_{k-1})$ is a square-free word over $\{\alpha,\beta,\gamma\}$.
  \end{itemize}
It is easy to verify that if there is a square in $\phi'$ then the colors inherited from $\phi$ form a square on $G$. 
\end{proof}

\begin{theorem}\label{mainth}
  Let $G$ be a graph and $c= \max\left\{35216,8\frac{\log(2\pi'(G)+3)}{\log(1.01)}+216\right\}$.
  Then for any $(\ge c)$-subdivision $H$ of $G$, 
  $$\pi(H)=3\,.$$  
\end{theorem}
\begin{proof}
Let $n=\max\left\{8750,2\frac{\log(2\pi'(G)+3)}{\log(1.01)}\right\}$. This implies $c= 4n+216$.

Let $G'$ be a subdivision of $G$ such that $H$ is a $(\ge 4n+216, \le 9n)$-subdivision of $G$'.
Let us first show that there exists such a graph $H$.
Since $4n+216<2\times (9n)$, for any integer $x\ge4n+216$ there exists an integer $\gamma(x)$ such that
$4n+216\le\frac{x}{\gamma(x)}\le 9n$.
Thus, for any edge $e$ of $G$ that is subdivided $k$ times in $H$, we can choose $e$ to be subdivided $\frac{k}{\gamma(k)}$ times in $G'$.

By Lemma \ref{lemmasubddoesnotincreasepi},  $2\pi'(G)+3\ge \pi'(G')$.
Let us now show that we can apply Lemma \ref{lemmaConstruction} to $G'$ and $H$.
We can verify by simple computation that $16(60n^2+1)< 1.01^{n/2}$ for any $n\ge8750$.
Hence, by definition of $n$.
$$\frac{1.01^{n}}{16(60n^2+1)}\ge1.01^{n/2} \ge 2\pi'(G)+3\ge \pi'(G')\,.$$
So $n$ verifies the conditions of Lemma \ref{lemmaConstruction}. Since $H$ is a 
$(\ge 4n+216, \le 9n)$-subdivision of $G'$ we can apply Lemma \ref{lemmaConstruction} and we  conclude that 
$$\pi(H)=3\,.\qedhere$$
\end{proof}

\section{Improving the coefficient}
We showed in Theorem \ref{mainth} that for any graph $G$ as long as there are at least $558 \log_2(\pi'(G))+ O(1)$ division vertices per edge the resulting graph is $3$-colorable.
This result is optimal in the sense that $\Omega(\log n)$ division vertices are needed  on some edges of any nonrepetitively $O(1)$-colorable subdivision of $K_n$\cite{nesetrilsubdivision}.
However, we can try to reduce the multiplicative constant $558$.
By being more careful on the computations, we can replace $558$ by 
$\frac{4}{\log_2(\gamma(\mathcal{N}))}$ where $\gamma(\mathcal{N})$ is the growth rate of the set of nice words. By using our lower bound of $1.01$ on the growth rate of the number of nice words (Theorem \ref{expboundsqpfree}) we obtain the bound $279\log_2(n)+ O(1)$. But we expect the growth rate to be much closer to $1.3$.

In fact, if instead of $\pal=\overline{\hpal} 0 \hpal$ we take any long enough square-free palindrome $\pal'=\overline{\hpal}' 0 \hpal'$ then it is easy to adapt the proof from \cite{shursqf} to show that the growth rate of the set of square-free words that avoids $\hpal'$ and $\overline{\hpal}'$ can be arbitrarily close to $1.3$. It probably does not change the growth rate to add the constraint that $\pal' w\pal'$ be square-free for every element $w$. However, we do not know how to prove this second point, but if it holds we can then replace $1.01$ by $1.3$. Our lower bound on the number of required division vertices per edge becomes $10.56 \log_2(n)+ O(1)$. We suspect that this coefficient would still be far from optimal.

\bibliographystyle{plain}

\begin{thebibliography}{1}

\bibitem{alongraph}
N.~Alon, J.~Grytczuk, M.~Haluszcza, and O.~Riordan.
\newblock Nonrepetitive colorings of graphs.
\newblock {\em Random Structures \& Algorithms}, 21:336--346, 2002.

\bibitem{3subdivision4colors}
J.~Barát and D.~R.~Wood.  
\newblock Notes on nonrepetitive graph colour-ing.
\newblock {\em Electronic Journal of Combinatorics}, 15:R99, 2008.

\bibitem{nicesubstitution}
J. Currie, T. Harju, P. Ochem, N. Rampersad,
\newblock Some further results on squarefree arithmetic progressions in infinite words,
\newblock {\em Theoretical Computer Science}, 799:140-148, 2019.


\bibitem{BerstelThue}
J.~Berstel.
\newblock Axel {T}hue's papers on repetitions in words: a translation.
\newblock {\em Publications du LaCIM 20, Universit\'e du Qu\'ebec à
  Montr\'eall}, 1995.

\bibitem{Dujmovic2016}
V.~Dujmovi{\'{c}}, G.~Joret, J.~Kozik, and D.~R. Wood.
\newblock Nonrepetitive colouring via entropy compression.
\newblock {\em Combinatorica}, 36(6):661--686, Dec 2016.

\bibitem{3uniformcoloring}
P.Erdős   and   L.  Lovász.
\newblock Problems  and  results  on 3-chromatic hypergraphs  and  some  related  questions.
\newblock  In {\em Infinite   and   Finite   Sets}, vol. 10 of {\em Colloq.  Math.  Soc.  János  Bolyai}, pp. 609–627. North-Holland, 1975

\bibitem{MontassierEntropie}
D.~Gon\c{c}alves, M.~Montassier, and A.~Pinlou.
\newblock Entropy compression method applied to graph colorings.
\newblock {\em arXiv e-prints}, arXiv:1406.4380, 2014.

\bibitem{Grytczuksubdivisionsurvey}
J. Grytczuk.   
\newblock Nonrepetitive colorings of graphs—a survey.
\newblock {\em Int.J. Math. Math. Sci.},  74639,  2007.

\bibitem{HARANT2012374}
J.~Harant and S.~Jendrol.
\newblock Nonrepetitive vertex colorings of graphs.
\newblock {\em Discrete Mathematics}, 312(2):374--380, 2012.

\bibitem{nesetrilsubdivision}
J.~Nešetřil, P.~Ossona de Mendez, and D.~R. Wood.
\newblock Characterisations and examples of graph classes with bounded expansion.
\newblock {\em Euro-pean J. Combin.}, 33(3):350–373, 2011.

\bibitem{3subdivision}
A. Pezarski and M. Zmarz.
\newblock Non-repetitive 3-coloring of subdivided graphs.
\newblock {\em Electronic Journal of Combinatorics}, 16(1), 2009.

\bibitem{rosenfeldCounting}
M.~Rosenfeld.
\newblock Another approach to non-repetitive colorings of graphs of bounded degree.
\newblock {\em Electronic Journal of Combinatorics}, 27(3), 2020.

\bibitem{shursqf}
A. M.~Shur.
\newblock Two-Sided Bounds for the Growth Rates of Power-Free Languages.
\newblock In: Diekert V., Nowotka D. (eds){\em Developments in Language Theory}, DLT 2009.
 Lecture Notes in Computer Science, vol 5583.
 
\bibitem{Thue06}
A.~Thue.
\newblock {\"{U}ber unendliche {Z}eichenreihen}.
\newblock {\em 'Norske Vid. Selsk. Skr. I. Mat. Nat. Kl. Christiania}, 7:1--22,
  1906.

\bibitem{Thue1}
A.~Thue.
\newblock {\"Uber} die gegenseitige {L}age gleicher {T}eile gewisser
  {Z}eichenreihen.
\newblock {\em Norske Vid. Selsk. Skr. I. Mat. Nat. Kl. Christiania,},
  10:1--67, 1912.

\bibitem{woodSurvey}
D.~R. Wood.
\newblock Nonrepetitive Graph Colouring.
\newblock {\em arXiv e-prints}, arXiv:2009.02001, 2020.
\end{thebibliography}

\end{document}